\documentclass[10pt,reqno]{amsart}

\usepackage[all]{xy}
\usepackage{a4wide}
\usepackage{amssymb}
\usepackage[mathscr]{eucal}

\newtheorem{prop}{Proposition}[section]
\newtheorem{lem}[prop]{Lemma}
\newtheorem{thm}[prop]{Theorem}
\newtheorem{cor}[prop]{Corollary}

\theoremstyle{remark}

\theoremstyle{definition}
\newtheorem{defn}[prop]{Definition}

\numberwithin{equation}{section}
\numberwithin{prop}{section}

\DeclareMathOperator{\Mor}{Mor}
\DeclareMathOperator{\Aut}{Aut}
\DeclareMathOperator{\M}{M}
\DeclareMathOperator{\C}{C}
\DeclareMathOperator{\Co}{C_\infty}
\DeclareMathOperator{\spec}{Sp}

\newcommand{\CC}{\mathbb{C}}
\newcommand{\RR}{\mathbb{R}}
\newcommand{\comp}{\!\circ\!}
\newcommand{\tens}{\otimes}
\newcommand{\cst}{\mathrm{C}^*}
\newcommand{\Del}{\Delta}
\newcommand{\id}{\mathrm{id}}

\newcommand{\ZZ}{{\mathbb Z}}
\newcommand{\TT}{{\mathbb T}}
\newcommand{\zz}{{\boldsymbol z}}
\newcommand{\zzb}{\overline{\zz}}
\newcommand{\bez}{\setminus}
\newcommand{\eps}{\varepsilon}
\newcommand{\ph}{\varphi}
\newcommand{\la}{\langle}
\newcommand{\ra}{\rangle}
\newcommand{\aff}{\,\,\!\eta\,\,\!}
\newcommand{\dplus}{\,\dot{+}\,}
\newcommand{\Y}{{\Gamma}}
\newcommand{\Ybar}{{\overline{\Y}}}
\newcommand{\Yhat}{\widehat{\Y}}
\newcommand{\y}{\gamma}
\newcommand{\yhat}{\widehat{\y}}

\newcommand{\zbar}{\overline{z}}
\newcommand{\alhat}{\widehat{\alpha}}
\newcommand{\behat}{\widehat{\beta}}
\newcommand{\ghat}{\hat{g}}
\newcommand{\DelG}{\Del_{\scriptscriptstyle G}}
\newcommand{\DelYhat}{\Del_{\scriptscriptstyle\Yhat}}
\newcommand{\PhiG}{\Phi_{\scriptscriptstyle G}}
\newcommand{\Hs}{G\bigl/\Yhat\bigr.}

\newcommand{\CYbar}{\C\bigl(\Ybar\bigr)}
\newcommand{\CoYbar}{\Co\bigl(\Ybar\bigr)}
\newcommand{\CoYhat}{\Co\bigl(\Yhat\bigr)}
\newcommand{\CoHs}{\Co\bigl(\Hs\bigr)}
\newcommand{\cpY}{\CoYbar\!\!\rtimes_{\beta}\!\Y}
\newcommand{\Int}{\int\limits}
\newcommand{\aA}{\mathscr{A}}
\newcommand{\bB}{\mathscr{B}}
\newcommand{\aAb}{\mathscr{A}_{\textrm{\tiny b}}}
\newcommand{\bBb}{\mathscr{B}_{\textrm{\tiny b}}}
\newcommand{\cont}{\cdot}
\newcommand{\kaptil}{\widetilde{\kappa}}
\newcommand{\Omtil}{\widetilde{\Omega}}
\newcommand{\om}{\omega}
\newcommand{\omb}{\overline{\omega}}
\newcommand{\DLz}[1]{
\displaystyle{\frac{\partial_{\scriptscriptstyle\mathrm{L}}{#1}}{\partial\zz}}}
\newcommand{\DLzb}[1]{
\displaystyle{\frac{\partial_{\scriptscriptstyle\mathrm{L}}{#1}}{\partial\zzb}}}
\newcommand{\DRz}[1]{
\displaystyle{\frac{\partial_{\scriptscriptstyle\mathrm{R}}{#1}}{\partial\zz}}}
\newcommand{\DRzb}[1]{
\displaystyle{\frac{\partial_{\scriptscriptstyle\mathrm{R}}{#1}}{\partial\zzb}}}
\newcommand{\Cinf}{\C^\infty_q\bigl(\Ybar\bigr)}
\newcommand{\LYbar}{L^2\bigl(\Ybar\bigr)}
\newcommand{\is}[2]{\left(#1\:\vline\:#2\right)}
\newcommand{\g}[1]{g_{\scriptscriptstyle{#1}}}
\newcommand{\F}{\mathscr{F}}
\newcommand{\Fq}{F_{\!\scriptscriptstyle{q}}}
\newcommand{\zZ}{\mathscr{Z}}
\newcommand{\DS}{\displaystyle}
\newcommand{\LY}{L^2(\Y)}
\newcommand{\LinfYbar}{L^\infty\bigl(\Ybar\bigr)}

\begin{document}

\title{Analysis on a homogeneous space of a quantum group}

\date{\today}

\author{W.~Pusz}
\address{Department of Mathematical Methods in Physics\\
Faculty of Physics\\
Warsaw University}
\email{wieslaw.pusz@fuw.edu.pl}

\author{Piotr M.~So{\l}tan}
\address{Department of Mathematical Methods in Physics\\
Faculty of Physics\\
Warsaw University}
\email{piotr.soltan@fuw.edu.pl}

\thanks{Research partially supported by KBN grants
nos.~115/E-343/SPB/6.PRUE/DIE50/2005-2008,
2P03A04022 \& 1PO3A03626}

\begin{abstract}
A detailed account of the construction of a homogeneous space for the quantum
``$az+b$'' group is presented. The homogeneous space is described by a
commutative $\mathrm{C}^*$-algebra which means that it is a classical space.
Then a
covariant differential calculus on the homogeneous space is constructed and
studied. A covariant measure and an analogue of the exponential function are
used to introduce elements o Fourier analysis.
\end{abstract}

\maketitle

\section{Introduction}

The concept of a homogeneous space for quantum groups has been studied for some
time now. In case of compact quantum groups the definition and many examples
were worked out by Podleś (\cite{podles}). For non compact quantum groups even
the definition is not agreed upon by the experts. In order to have a better
understanding of the problems we introduce an example of an object which
should, in our opinion, be called a homogeneous space. It turns out that
despite a fairly abstract approach to the definition of our homogeneous space
for the quantum ``$az+b$'' group for real deformation parameter
(cf.~\cite{azb}), the object turns out to be a classical space, i.e.~is is a
quantum space described by a commutative $\cst$-algebra. This relatively simple
situation makes it easier to deal with such aspects of quantum group
covariant non commutative geometry as differential calculus, covariant
measures. Using some tools previously employed for the construction of examples
of non compact quantum groups we are able to introduce elements of Fourier
analysis and prove that our ``quantum Fourier transform'' has properties
similar to those of its classical counterpart.

The paper is organizes as follows. In Section \ref{AzB} we recall the necessary
information about the quantum ``$az+b$'' group needed for our construction. The
homogeneous is defined and described in detail in Section \ref{ths}. Then we
introduce a covariant differential calculus on the homogeneous space in Section
\ref{diffCal}. At first we describe a general construction of, so called,
embeddable covariant bimodules for Hopf algebras. Then we introduce the class
of function on the homogeneous space which play the role of smooth functions in
classical setting. In Subsection \ref{calcu} we give an abstract description of
our covariant differential calculus and in the next subsection we describe its
concrete realization and covariance properties. Section \ref{covMeas} is
devoted to developing integration on the homogeneous space. A measure is
introduced and properties of differential operators coming from the
differential calculus are studied from the point of view of functional
analysis. The covariance and uniqueness of the chosen measure are described in
Subsection \ref{covMu}. Finally Section \ref{FTs} contains the definition and
study of basic properties of the analogue of Fourier transform on the
homogeneous space.

\section{Quantum ``$az+b$'' group}\label{AzB}

Let us briefly recall the basic facts about the quantum ``$az+b$'' group
(cf.~\cite[Appendix A]{azb}). For a real parameter $q$ such that
$0<q<1$ consider
\[
\Y=\left\{z\in\CC:\:|z|\in q^{\ZZ}\right\}.
\]
Clearly $\Y$ is a multiplicative subgroup of $\CC\bez\{0\}$ and
$\Y\simeq\ZZ\times\TT$. Moreover $\Y$ is self dual, $\Y\simeq\Yhat$.
Any $\y\in\Y$ is of the form $\y=q^{i\ph+k}$ for unique $k\in\ZZ$ and
$\ph\in\left[0,-\frac{2\pi}{\log{q}}\right[$. Let
\[
\chi(\y,\y')=\chi\left(q^{i\ph+k},q^{i\ph'+k'}\right)=q^{i(\ph k'+\ph'k)}
\]
for all $\y,\y'\in\Y$. Then $\chi:\Y\times\Y\to\TT$ is a non degenerate
bicharacter on $\Y$. Using the non degeneracy of $\chi$ we shall identify
$\Yhat$ with $\Y$ via the formula
\[
\la\yhat,\y\ra=\chi(\yhat,\y).
\]
Finally let $\Ybar$ be the closure of $\Y$, $\Ybar=\Y\cup\{0\}$.

The quantum ``$az+b$'' group is a pair $G=(A,\DelG)$, where $A$ is a
$\cst$-algebra corresponding to the algebra of all continuous functions
vanishing at infinity on the group and comultiplication $\DelG$ is a
coassociative morphism: $\DelG\in\Mor(A,A\tens A)$, i.e.~a $*$-homomorphism
from $A$ to $\M(A\tens A)$ encoding the group structure.

According to \cite{azb} there are normal elements $a$ and $b$ affiliated
with the $\cst$-algebra $A$ such that $\spec{a},\:\spec{b}\subset\Ybar$.
Furthermore $a$ is invertible and $a^{-1}\aff A$. Therefore for all $\y\in\Y$
we can form unitary elements $\chi(a,\y)$ of $\M(A)$ using functional calculus.
Moreover $a$ and $b$ satisfy the relation
\begin{equation}\label{weyl}
\chi(a,\y)b\chi(a,\y)^*=\y b
\end{equation}
for any $\y\in\Y$. It turns out that $A$ is generated by $a$, $a^{-1}$ and
$b$ (cf.~\cite{gen}). Let us also note that formula \eqref{weyl} implies that
$ab$, $ba$, $ab^*$ and $b^*a$ are well defined normal operators satisfying
\begin{equation}\label{algrel}
ab=q^2ba,\quad ab^*=b^*a
\end{equation}
(cf.~\cite{azb}). The reader should be warned, however, that \eqref{weyl} and
\eqref{algrel} are not equivalent.

Let us remark (cf.~\cite{gen}) that $A$ is the universal $\cst$-algebra
generated by $a$ and $b$ satisfying the described relations in the sense that
for any $\cst$-algebra $C$ and $a_0,b_0\aff C$ satisfying
\begin{equation}\label{ComRel}
\left.\begin{array}{l}
a_0,b_0\textrm{ are normal},\\
\spec{a_0},\spec{b_0}\subset\Ybar,\\
a_0\textrm{ is invertible and }a_0^{-1}\aff C\\
\chi(a_0,\y)b_0\chi(a_0,\y)^*=\y b_0\textrm{ for all }\y\in\Y
\end{array}\right\}
\end{equation}
there exists a unique morphism $\Psi\in\Mor(A,C)$ such that
\[
a_0=\Psi(a),\quad b_0=\Psi(b).
\]

It turns out that the $\cst$-algebra $A$ has a relatively simple structure. It
is a $\cst$-crossed product
\begin{equation}\label{crossed}
A=\cpY,
\end{equation}
where the action
$\beta:\Y\ni\y\mapsto\beta_{\y}\in\Aut\bigl(\CoYbar\bigr)$ is given by
\begin{equation}\label{dzial}
\left(\beta_{\y}f\right)(\y')=f(\y'\y)
\end{equation}
for all $f\in\CoYbar$.

To describe the generating elements $a$ and $b$ in this setting note that
the natural inclusion $\CoYbar\hookrightarrow\M\bigl(\cpY\bigr)$ is a
morphism form $\CoYbar$ to $\cpY$. Let $\zz$ be the standard
generator of $\CoYbar$, i.e.~$\zz(\y)=\y$ for all $\y\in\Ybar$ then
$b\aff\cpY$ is the image of $\zz$ under the canonical inclusion.

By definition of a crossed product $\M\bigl(\cpY\bigr)$ contains a strictly
continuous family of unitaries $\left(U_{\y}\right)_{\y\in\Y}$ such that
\[
\beta_\y(f)=U_{\y}fU_{\y}^*
\]
for any $f\in\CoYbar$. The results of Sections 4.~and 5.~of \cite{azb} show
that there exists a normal element $a\aff\cpY$ such that
$\spec{a}\subset\Ybar$, $\ker{a}=\{0\}$ and
\begin{equation}\label{Uy}
U_\y=\chi(a,\y)
\end{equation}
for all $\y\in\Y$. Moreover $a^{-1}\aff\cpY$. The elements
$a$, $a^{-1}$ and $b$ generate $\cpY$ in the sense of \cite{gen}.

The comultiplication is defined in the following way. Consider the $\cst$-tensor
product $A\tens A$. It turns out that $a\tens b+b\tens I$ is a closable
operator and its closure $a\tens b\dplus b\tens I$ is a normal element
affiliated with $A\tens A$. Moreover
$\spec(a\tens b\dplus b\tens a)\subset\Ybar$. Clearly
$(a\tens a)\aff A\tens A$ is normal, invertible,
$(a^{-1}\tens a^{-1})\aff A\tens A$ and $\spec{(a\tens a)}\subset\Ybar$. One
can check that $a_0=a\tens a$ and $b_0=a\tens b\dplus b\tens I$ satisfy the
commutation relations \eqref{ComRel}. Therefore there exists a unique
$\DelG\in\Mor(A,A\tens A)$ such that
\begin{equation}\label{delG}
\begin{array}{r@{\;=\;}l}
\DelG(a)&a\tens a,\\
\DelG(b)&a\tens b\dplus b\tens I.
\end{array}
\end{equation}
Moreover $\DelG$ is coassociative. This completes the description of the
quantum ``$az+b$'' group on $\cst$-algebra level.

\section{The homegeneous space}\label{ths}

In this section we shall introduce a homegeneous space for $G$ which is the
main object of our considerations. First we shall observe that $\Yhat$ is a
subgroup of $G$. Indeed: for $\yhat\in\Yhat$ let
\begin{equation}\label{wartpi}
\begin{array}{r@{\;=\;}l}
a_0(\yhat)&\yhat,\\
b_0(\yhat)&0.
\end{array}
\end{equation}
Then $a_0$ and $b_0$ are continuous functions on $\Yhat$, i.e.~elements
affiliated with the $\cst$-algebra $\CoYhat$. They satisfy the relations
\eqref{ComRel} and therefore there exists a unique morphism
$\pi\in\Mor\bigl(A,\CoYhat\bigr)$ such that $a_0=\pi(a)$ and $b_0=\pi(b)$.
Clearly $\pi$ is surjective and the reader will easily check that
\begin{equation}\label{subgrp}
(\pi\tens\pi)\comp\DelG=\DelYhat\comp\pi,
\end{equation}
where $\DelYhat$ is the standard comultiplication on $\CoYhat$,
i.e.~$\left(\DelYhat f\right)(\yhat_1,\yhat_2)=f(\yhat_1\yhat_2)$. This means
that $\Yhat$ is a subgroup of $G$.

Our aim in this section is to describe a quantum analogue of the homogeneous
space $\Hs$. In the classical situation when $G$ is a locally compact group
and $\Yhat$ is a subgroup of $G$, continuous functions on $\Hs$ may be
identified with continuous functions on $G$ which are constant on left cosets,
i.e.~such $x\in\C{G}$ that
\begin{equation}\label{stale}
x(g\yhat)=x(g)
\end{equation}
for all $g\in G$ and $\yhat\in\Yhat$. This space of functions carries a natural
left action of $G$ by left shifts. In contrast to the case of compact quantum
groups (\cite[Section 1]{podles}) at the moment there seems to be no
appropriate definition of a homogeneous space for non compact quantum groups in
the $\cst$-algebra approach. The main problem is to describe the class of
``functions on $G$'' corresponding to continuous functions vanishing at
infinity on $\Hs$. Clearly they should be bounded and continuous.

In the case of the classical ``$az+b$'' group $G$ is topologically the cartesian
product $G=\Yhat\times\left(\Hs\right)$, where $\Yhat$ is the subgroup of
homoteties ($b=0$). Therefore in this case a bounded continuous function $x$ on
$G$ corresponds to a function vanishing at infinity on $\Hs$ if and only
if $x$ satisfies \eqref{stale} and the function
\begin{equation}\label{znikanie}
G\ni g=(a,b)\longmapsto x(a,b)f(a)\in\CC
\end{equation}
vanishes at infinity on $G$ for any $f\in\CoYhat$. In other words $x(a,b)f(a)$
belongs to $\Co(G)$.

We shall follow the above ideas in the case of the {\em quantum}\/ ``$az+b$''
group.
Any ``continuous function'' on $\Hs$ is realized by an element $x\aff A$ such
that (cf.~\eqref{stale})
\begin{equation}\label{Stale}
(\id\tens\pi)\DelG(x)=x\tens I.
\end{equation}
One can check that $x=b$ is a solution of Equation \eqref{Stale}.
Therefore for any $f\in\CYbar$ the element $f(b)$
affiliated with $A$ satisfies \eqref{Stale}. Since the $\cst$-algebra $A$ is
generated by $a$, $a^{-1}$ and $b$, and $x=a$ does not fulfill the
requirement \eqref{Stale}, one expects that the algebra of all continuous
functions vanishing at infinity on $\Hs$ coincides with
$\bigl\{f(b):\:f\in\CoYbar\bigr\}$ and the quotient map is given by the
morphism
\begin{equation}\label{quotient}
\CoYbar\ni f\longmapsto f(b)\in\M(A).
\end{equation}

\begin{defn}\label{HS}
We define $\CoHs$ to be the set of those $x\in\M(A)$ which satisfy
\begin{enumerate}
\item\label{l1} $(\id\tens\pi)\DelG(x)=x\tens I$,
\item\label{l2} $xf(a)\in A$ for all $f\in\CoYhat$

\noindent and
\item\label{l3} the map $\Y\ni\y\mapsto U_{\y}xU_{\y}^*\in\M(A)$ is
norm continuous.
\end{enumerate}
\end{defn}

We shall refer to \eqref{l3} of Definition \ref{HS} as the {\em regularity
condition.} This condition was not apparent in the classical setting and is of
purely quantum nature. In introducing this condition we followed the idea of
Landstad (\cite{lan}).

\begin{prop}
$\CoHs$ is a $\cst$-algebra.
\end{prop}

\begin{proof}
It is clear that the set of elements $x\in\M(A)$ satisfying \eqref{l1} and
\eqref{l3} of Definition \ref{HS} is a $\cst$-algebra. Thus it is enough to
show that if $x$ satisfies \eqref{l2} and \eqref{l3} of Definition \ref{HS}
then so does $x^*$.

It follows from Lemma \ref{LemSLW} (see below) that for any $f\in\CoYhat$ the
element
\[
x^*f(a)=\bigl(\overline{f}(a)x\bigr)^*
\]
belongs to $A$. Therefore $\CoHs$ is a $\cst$-algebra.
\end{proof}

For any $g\in L^1(\Y)$ let $\ghat$ denote its Fourier transform:
\begin{equation}\label{F}
\ghat(\yhat)=\Int_{\Y}g(\y)\la\yhat,\y\ra\,d\y
=\Int_{\Y}g(\y)\chi(\yhat,\y)\,d\y,
\end{equation}
where $d\y$ is a fixed Haar measure on $\Y$. With this notation we have
\begin{equation}\label{FF}
\ghat(a)=\Int_{\Y}g(\y)\chi(a,\y)\,d\y=\Int_{\Y}g(\y)U_{\y}\,d\y.
\end{equation}
Clearly $\ghat(a)\in\M(A)$.

\begin{lem}\label{LemSLW}
Let $x\in\M(A)$ be such that the map $\Y\ni\y\mapsto U_{\y}xU_{\y}^*$ is
norm continuous and assume that $xf(a)\in A$ for all $f\in\CoYhat$. Then
$f(a)x\in A$ for all $f\in\CoYhat$.
\end{lem}

\begin{proof}
By the continuity of the map $\Y\ni\y\mapsto U_{\y}xU_{\y}^*$, for any
$\eps>0$ there exists a compact neighborhood $V_{\eps}$ of $1\in\Y$ such that
\[
\|U_{\y}x-xU_{\y}\|<\eps
\]
for any $\y\in V_{\eps}$. We shall choose $V_\eps$ in such a way that
$V_{\eps}$ shrinks to $\{1\}$ when $\eps\to 0$.
Let
\[
g_{\eps}(\y)=\left\{\begin{array}{@{\:}c@{\quad\textrm{when}\;}l}
\bigl(\textrm{Haar measure of }V_\eps\bigr)^{-1}&\y\in V_\eps,\\
0&\y\not\in V_\eps.\end{array}\right.
\]
and let $\ghat_{\eps}$ be the Fourier transform of $g_\eps$ (cf.~\eqref{F}):
\[
\ghat_{\eps}(\yhat)=\Int_{\Y}g_{\eps}(\y)\chi(\yhat,\y)\,d\y.
\]
Using \eqref{FF} with $g_\eps$ instead of $g$ we obtain
\[
\|\ghat_{\eps}(a)x-x\ghat_{\eps}(a)\|=
\left\|\Int_{\Y}g_{\eps}(\y)(U_{\y}x-xU_{\y})\,d\y\right\|
\leq\Int_{\Y}g_{\eps}(\y)\left\|U_{\y}x-xU_{\y}\right\|\,d\y\leq\eps.
\]
Therefore for any $f\in\CoYhat$ we have
\[
\|f(a)\ghat_{\eps}(a)x-f(a)x\ghat_{\eps}(a)\|\leq\|f(a)\|\eps.
\]
By construction $f(a)\ghat_{\eps}(a)\to f(a)$ in norm when
$\eps\to0$.

Since
\[
\left\|f(a)x-f(a)x\ghat_{\eps}(a)\right\|\leq
\left\|\bigl(f(a)-f(a)\ghat_{\eps}(a)\bigl)x\right\|
+\left\|f(a)\ghat_{\eps}(a)x-f(a)x\ghat_{\eps}(a)\right\|
\xrightarrow[\eps\to0]{}0
\]
and $f(a)x\ghat_{\eps}(a)=f(a)\bigl(x\ghat_{\eps}(a)\bigr)\in\M(A)A\subset A$
we see that $f(a)x\in A$.
\end{proof}

Now we have the following

\begin{thm}\label{HomSpace}
$\CoHs=\bigl\{f(b):\:f\in\CoYbar\bigr\}$.
\end{thm}

\begin{proof}
For $\yhat\in\Yhat$ let $\pi_{\yhat}$ be the morphism
$\pi\in\Mor\bigl(A,\CoYhat\bigr)$
(introduced in the beginning of this Section) composed with evaluation at
$\yhat$, $\pi_{\yhat}(x)=\left(\pi(x)\right)(\yhat)$. Then
$\pi_{\yhat}\in\Mor(A,\CC)$ and due to \eqref{subgrp}
\begin{equation}\label{alhat2}
\pi_{\yhat_1}*\pi_{\yhat_2}=(\pi_{\yhat_1}\tens\pi_{\yhat_2})\comp\DelG
=\pi_{\yhat_1\yhat_2}.
\end{equation}
For $x\in A$ and $\yhat\in\Yhat$ let
\begin{equation}\label{alhat}
\alhat_{\yhat}(x)=\pi_{\yhat}*x=(\id\tens\pi_{\yhat})\DelG(x).
\end{equation}
Then by \eqref{alhat2}, $\alhat_{\yhat}$ is an automorphism of $A$ and
$\Yhat\ni\yhat\mapsto\alhat_{\yhat}\in\Aut{A}$ is a continuous action of $\Yhat$
on $A$. In other words $\left(A,\Yhat,\alhat\right)$ is a
$\cst$-dynamical system. Clearly (cf.~\eqref{alhat}) for any $x\aff A$
\begin{equation}\label{Equiv}
\Bigl((\id\tens\pi)\DelG(x)=x\tens I\Bigr)\Longleftrightarrow
\left(\begin{array}{c}\alhat_{\yhat}(x)=x\text{ for}\\
\text{all }\yhat\in\Yhat\end{array}\right)
\end{equation}
and since by \eqref{delG}
\[
\DelG\left(U_\y\right)=\DelG\left(\chi(a,\y)\right)=\chi(a\tens a,\y)=
\chi(a,\y)\tens\chi(a,\y)=U_\y\tens U_\y
\]
we have
\[
\alhat_{\yhat}\left(U_{\y}\right)=\la\yhat,\y\ra U_{\y}
\]
for all $\y\in\Y$ and $\yhat\in\Yhat$.

\newlength{\nl}
\settowidth{\nl}{\qquad{\small (\emph{iii}) the map
$\Y\ni\y\mapsto U_{\y}xU_{\y}^*\in\M(A)$is norm continuous}}

By Landstad's theorem (see e.g.~\cite[Theorem 7.8.8]{ped}) there exists a
$\cst$-dynamical system $(B,\Y,\alpha)$ such that $A=B\rtimes_{\alpha}\Y$.
Moreover
\begin{equation}\label{LC}
B=\left\{x\in\M(A):
\begin{minipage}[c]{\nl}{\small
\begin{enumerate}
\item[(\emph{i})] $\alhat_{\yhat}(x)=x$ for all $\yhat\in\Yhat$,
\item[(\emph{ii})] $x\ghat(a),\ghat(a)x\in A$ for all $g\in L^1(\Y)$,
\item[(\emph{iii})] the map $\Y\ni\y\mapsto U_{\y}xU_{\y}^*\in\M(A)$
is norm continuous
\end{enumerate}}
\end{minipage}
\right\}.
\end{equation}
The action $\alpha$ of $\Y$ on $B$ is given by
\[
\alpha_{\y}(x)=U_{\y}xU_{\y}^*.
\]
Furthermore $(B,\Y,\alpha)$ is unique up to covariant isomorphism and thus by
\eqref{crossed} we have $(B,\Y,\alpha)\simeq\bigl(\CoYbar,\Y,\beta\bigr)$ . It
turns out that in our case
\[
B=\CoHs
\]
as subsets of $\M(A)$. The condition (\emph{i}) from \eqref{LC} coincides with
\eqref{l1} of Definition \ref{HS} due to \eqref{Equiv}. Condition
\eqref{l3} of Definition \ref{HS} is the same as (\emph{iii}) from \eqref{LC}.
At first sight (\emph{ii}) from \eqref{LC} seems stronger than \eqref{l2} of
Definition \ref{HS}, but in fact they are equivalent by Lemma \ref{LemSLW}.

Now let $x=f(b)$ with $f\in\CoYbar$. Clearly $x\in\M(A)$ and $x$
satisfies ({\em i}\/). Since
\begin{equation}\label{cov}
U_{\y}f(b)U_{\y}^*=f\left(U_{\y}bU_{\y}^*\right)
=f(\y b)=\left(\beta_{\y}(f)\right)(b)
\end{equation}
condition ({\em iii}\/) of \eqref{LC} is satisfied.

It is known (cf.~\cite[Formula (4.3)]{azb}) that $f(b)g(a)\in A$ for all
$f\in\CoYbar$ and $g\in\CoYhat$. In particular $xg(a)\in A$ for all
$g\in\CoYhat$. By Lemma \ref{LemSLW} also $g(a)x\in A$ for any $f\in\CoYhat$.
It means that $x=f(b)$ satisfies condition ({\em ii}\/) of \eqref{LC}.

As a result $\bigl\{f(b):\:f\in\CoYbar\bigr\}\subset B$.

Note that due to \eqref{cov} the inclusion map
\[
\CoYbar\ni f\longmapsto f(b)\in B
\]
intertwines the actions $\beta$ and $\alpha$ on $\CoYbar$ and $B$. By
composing this map with the covariant isomorphism $B\simeq\CoYbar$ we obtain
a covariant inclusion $\jmath:\CoYbar\hookrightarrow\CoYbar$. Therefore the
underlying map
\[
\Psi:\Ybar\longrightarrow\Ybar.
\]
is surjective and equivariant for the natural action of $\Y$ on
$\Ybar$.\footnote{
The reader will notice that $\jmath$ is not a morphism of $\cst$-algebras,
so the existence of the map $\Psi$ is not obvious. However any character of $B$
can be treated as an irreducible representation of $\jmath(B)$ seen as a
subalgebra of $B$. This representation can be extended to $B$ as
an irreducible representation (\cite[Proposition 2.10.2]{dix}) which corresponds
to a point of $\Ybar$. Thus the map $\Psi$ exists and is surjective. It is,
however, a map from $\Ybar$ to its one point compactification and it is not
guaranteed to be continuous. Nevertheless it clearly has all of $\Ybar$ in its
range.}

Now for any $\y\in\Ybar$
\[
\Psi(\y)=\Psi(1\cdot\y)=\Psi(1)\y
\]
and $\Psi(1)\neq0$ because $\Psi$ is surjective. This shows that $\Psi$
is a homeomorphism and consequently $B=\bigl\{f(b):\:f\in\CoYbar\bigr\}$.
\end{proof}

Let us remark that Theorem \ref{HomSpace}
proves that the quantum homogeneous space $\Hs$ is in fact a classical
space: $\Hs=\Ybar$. More precisely $\CoHs=\CoYbar$
and the morphism $\CoYbar\ni f\mapsto f(b)\in\M(A)$ (cf.~\eqref{quotient})
is dual to the quotient map.

Using this morphism we shall identify the generator $\zz$ of $\CoYbar$ with
$b\aff A$. Now we can define the action $\PhiG$ of $G=(A,\DelG)$ on the
homogeneous space $\Hs=\Ybar$ by restricting $\DelG$ to the image of
$\CoYbar$ in $\M(A)$. With such a definition we have (cf.~\eqref{delG}):
\[
\PhiG(\zz)=a\tens\zz\dplus b\tens I.
\]
Indeed using results of \cite[Formula (2.6)]{gen} one can show that
$\left(a\tens\zz\dplus b\tens I\right)\aff A\tens\CoYbar$.
It is normal and $\spec(a\tens\zz\dplus b\tens I)\subset\Ybar$.
Therefore $\PhiG\in\Mor\bigl(\CoYbar,A\tens\CoYbar\bigr)$. Moreover the
coassociativity of $\DelG$ leads to the following commutative diagram
\[
\xymatrix{
\CoYbar\ar[rr]^-{\PhiG}\ar[dd]_{\PhiG}&&
A\tens\CoYbar\ar[dd]^{\DelG\tens\id}\\ \\
A\tens\CoYbar\ar[rr]^-{\id\tens\PhiG}&&A\tens A\tens\CoYbar}
\]
which means that $G$ is acting on $\Ybar$. This justifies our definition of
the homogeneous space.

The homogeneous space $\Hs$ shall now play a role analogous to the role of the
complex plane acted upon by the classical ``$az+b$'' group. The fact that
$\Ybar$ is a closure in $\CC$ of a multiplicative subgroup of $\CC\bez\{0\}$
will not be important for what follows. We shall therefore from now on denote a
generic element of $\Ybar$ by the symbol $z$.

\section{Differential calculus}\label{diffCal}

\subsection{A general construction}\label{algdiff}

The following proposition describes a general method of constructing covariant
differential calculi on coideals of Hopf algebras (so called {\em embeddable
homogeneous spaces}\/ cf.~\cite[Definition 1.8]{podles}).

\begin{prop}\label{HopfProp}
Let $\left(\aA_0,\Del_0\right)$ be a Hopf algebra with counit $e$ and let
$\bB_0$ be a left coideal of $\aA_0$.
\begin{enumerate}
\item Let $\ph$ be a linear functional on $\bB_0$. Define
$d:\bB_0\to\aA_0$ by $d(b)=\ph*b$ and let $\Omega$ be the image of $d$.
Then
\begin{enumerate}
\item\label{df11} $\Omega$ is a left coideal of $\aA_0$ and the diagram
\begin{equation}\label{HopfDiag}
\xymatrix{
\bB_0\ar[rr]^-{d}\ar[dd]_{\Del_0}&&
\Omega\ar[dd]^{\Del_0}\\ \\
\aA_0\tens\bB_0\ar[rr]^-{\id\tens d}&&\aA_0\tens\Omega}
\end{equation}
is commutative;
\item if $\ph$ satisfies
\begin{equation}\label{counit}
\ph(bc)=\ph(b)e(c)+e(b)\ph(c)
\end{equation}
for all $b,c\in\bB_0$ then $(\Omega,d)$ is a first order differential calculus
over $\bB_0$ which by statement \eqref{df11} is left covariant under the action
of $(\aA_0,\Del_0)$.
\end{enumerate}
\item Let $d:\bB_0\to\aA_0$ be a linear map. Let $\Omega$ be the
image of $d$ and define $\ph=e\comp d$. Then
\begin{enumerate}
\item if $\Omega$ is a left coideal of $\aA_0$ and the diagram \eqref{HopfDiag}
is commutative then for any $b\in\bB_0$ we have $d(b)=\ph*b$;
\item if $(\Omega,d)$ is a first order differential calculus over $\bB_0$ then
$\ph$ satisfies \eqref{counit} for all $b,c\in\bB_0$.
\end{enumerate}
\end{enumerate}
\end{prop}

We omit the proof which is purely computational. The proposition has an obvious
analogue for $*$-calculi over $*$-algebras covariant under action of a Hopf
$*$-algebra.

In our study to the quantum homogeneous space of the quantum ``$az+b$''
group constructed in Section \ref{ths} we shall construct a differential
calculus which corresponds to the algebraic construction given in Proposition
\ref{HopfProp}.

\subsection{Smooth functions on $\Hs$}\label{smooth}

Let $\bB$ be the algebra of all continuous functions on $\Ybar$ which can be
continued analitically from any circle $\left\{z\in\Ybar:\:|z|=q^n\right\}$ to a
holomorphic function on $\CC\bez\{0\}$. If $f\in\bB$ then for any $z\in\Ybar$
the function
\[
\RR\ni t\longmapsto f\left(q^{it}z\right)\in\CC
\]
has holomorphic continuation to an entire function. We shall denote the value of
this continuation at $t=-i$ by $f(q\cont z)$. The algebra of bounded functions
contained in $\bB$ will be denoted by $\bBb$.

$\bB$ will play the role of smooth functions on the quotient space $\Hs=\Ybar$.
The symbol $\bB$ is chosen in analogy with the algebraic situation described in
Subsection \ref{algdiff}. We also need to define the algebra of smooth
functions on the quantum group $G$ which acts on $\Ybar$.

Let $\aA$ be the set of those elements affiliated with $A$ which are entire
analytic for the scaling group of $G$ (for the natural topology, cf.
\cite{gen}). Let $\aAb$ be
the intersection of $\aA$ with $\M(A)$.

We shall now describe the analogy between the algebraic context of Subsection
\ref{algdiff} and our situation. The algebra $\bB$ can be embedded in the set
of elements affiliated with $A$ via the map $\jmath:f\mapsto f(b)$. The
following algebraic statements:
\begin{enumerate}
\item $\aA_0$ is a Hopf $*$-algebra, in particular it has an antipode and for
each $x\in\aA_0$ we have $\Del_0(x)\in\aA_0\tens\aA_0$;
\item $\bB_0$ is a left coideal in $\aA_0$, i.e.~$\bB_0$ is a subset of $\aA_)$
and for each $x\in\bB_0$ the element $\Del_0(x)$ belongs to $\aA_0\tens\bB_0$;
\end{enumerate}
have the following counterparts in our situation:
\begin{enumerate}
\item\label{st1} $\aAb$ is a unital $*$-subalgebra of $\M(A)$ contained in the
domain of $\kaptil$ and stable under $\kaptil$, for each $x\in\aA$ the element
$\Del(x)$ is entire analytic for $\left(\tau_t\tens\tau_t\right)_{t\in\RR}$;
\item\label{st2} $\bB$ is a subset of $\aA$ and for each $x\in\bB$ the element
$\Del(x)$ is affiliated with $A\tens\CoYbar$ and is entire analytic for
$\left(\tau_t\tens\tau_t\right)_{t\in\RR}$;
\end{enumerate}

Let us comment on statements \eqref{st1} and \eqref{st2}. The first one is a
consequence of the properties of analytic generators of one parameter groups of
automorphisms and the fact that for each $t\in\RR$ we have
$(\tau_t\tens\tau_t)\comp\Del=\Del\comp\tau_t$. As for statement \eqref{st2},
the fact that $\bB\subset\aA$ follows from the fact that $\bB$ is the set of
elements which are entire analytic for the group
$\left(\beta_{q^{it}}\right)_{t\in\RR}$ and for any $t\in\RR$ we have
\[
\tau_t\comp\jmath=\jmath\comp\beta_{q^{2it}}.
\]

The differential calculus over $\bB$ we want to construct will be analogous to
one described in Subsection \ref{algdiff} with $\ph$ corresponding to the map
which sends $b$ and $b^*$ to $1$ and all higher powers of $b$ and $b^*$ to $0$.

\subsection{The calculus}\label{calcu}

Let $\Omtil$ be the $\bB$ bimodule generated by two elements $\om$ and $\omb$
with the defining relations\footnote{
To keep our notation less complicated we will sometimes be writing $f(z)$
instead of $f$.}
\[
\begin{array}{r@{\;=\;}l@{\smallskip}}
\om f(z)&f(q\cont qz)\om,\\
\omb f(z)&f(q\cont q^{-1}z)\omb
\end{array}
\]
for all $f\in\bB$. Clearly by putting $\om^*=\omb$ we obtain a
$\bB$-$*$-bimodule structure on $\Omtil$. Now we define the operator
$d:\bB\to\Omtil$ by
\[
df=\left[\om-\omb,f\right]=(\om-\omb)f-f(\om-\omb).
\]
Clearly $d$ satisfies the Leibniz identity. The identity function
$\zz:\Ybar\ni z\mapsto z\in\CC$ and its adjoint $\zzb$ belong to $\bB$ and
we can compute their differentials:
\begin{equation}\label{dzom}
\begin{array}{r@{\;=\;}l@{\smallskip}}
d\zz&(q^2-1)\zz\om,\\
d\zzb&(q^2-1)\omb\zzb.
\end{array}
\end{equation}
Now expressing $df$ as a combination of $d\zz$ and $d\zzb$ we get the following
expressions
\[
\begin{array}{r@{\;=\;}l@{\medskip}}
df&\DLz{f}d\zz+d\zzb\DRzb{f}\\
&d\zz\DRz{f}+d\zzb\DRzb{f}\\
&d\zz\DRz{f}+\DLzb{f}d\zzb\\
&\DLz{f}d\zz+\DLzb{f}d\zzb,
\end{array}
\]
where
\begin{equation}\label{diffop}
\begin{array}{r@{\;=\;}l@{\medskip}}
\DRz{f}&\DS\frac{f(q^{-1}\cont q^{-1}z)-f(z)}{(q^{-2}-1)z},\\
\DLz{f}&\DS\frac{f(z)-f(q\cont qz)}{(1-q^2)z},\\
\DLzb{f}&\DS\frac{f(q\cont q^{-1}z)-f(z)}{(q^{-2}-1)\zbar},\\
\DRzb{f}&\DS\frac{f(z)-f(q^{-1}\cont qz)}{(1-q^2)\zbar}.
\end{array}
\end{equation}
The above computations make sense for $z$ away from $0$. The space $\Cinf$ of
smooth functions on $\Ybar$ is defined as
\[
\Cinf=\left\{f\in\bB:\:\exists\:\lim\limits_{\Y\ni z\to0}Df(z)\text{ for
any product $D$ of operators \eqref{diffop}}\right\}.
\]
For any $f\in\Cinf$ the right hand sides of \eqref{diffop} define elements of
$\Cinf$.

The properties of analytic generators of one parameter groups show that for any
$f\in\Cinf$ we have
\begin{equation}\label{LR}
\begin{array}{r@{\;=\;}l@{\medskip}}
\DRz{f}(z)&\DLz{f}(q^{-1}\cont q^{-1}z),\\
\DRzb{f}(z)&\DLzb{f}(q^{-1}\cont qz)\\
\end{array}
\end{equation}
as well as
\begin{equation}\label{star}
\DLz{f^*}=\left(\DRzb{f}\right)^*\quad\text{and}
\quad\DLzb{f^*}=\left(\DRz{f}\right)^*,
\end{equation}
where $*$ denotes the standard involution on $\bB$ (complex conjugation).

Using the Leibniz rule for $d$ we derive the value of our differential
operators on products of functions. Indeed, if $f,g\in\Cinf$ then we have
\[
\begin{array}{r@{\;=\;}l@{\medskip}}
d(fg)&df\,g+f\,dg\\
&\left(\DRz{f}d\zz+d\zzb\DLzb{f}\right)g
+f\left(\DRz{g}d\zz+d\zzb\DLzb{g}\right)\\
&\left(\DRz{f}g(q^{-1}\cont q^{-1}z)+f\DRz{g}\right)d\zz
+d\zzb\left(\DLzb{f}g+f(q\cont q^{-1}z)\DLzb{g}\right)
\end{array}
\]
It follows that
\begin{equation}\label{poz}
\begin{array}{r@{\;=\;}l@{\medskip}}
\DRz{fg}&\DRz{f}g(q^{-1}\cont q^{-1}z)+f\DRz{g},\\
\DLzb{fg}&\DLzb{f}g+f(q\cont q^{-1}z)\DLzb{g}.
\end{array}
\end{equation}
Similarly
\begin{equation}\label{pozb}
\begin{array}{r@{\;=\;}l@{\medskip}}
\DLz{fg}&\DLz{f}g+f(q\cont qz)\DLz{g},\\
\DRzb{fg}&\DRzb{f}g(q^{-1}\cont qz)+f\DRzb{g}.
\end{array}
\end{equation}
From \eqref{star}, \eqref{poz} and \eqref{pozb} we get

\begin{prop}
$\Cinf$ is a unital $*$-algebra.
\end{prop}

Define
\[
\Omega^1=\left\{\alpha d\zz\alpha'+\beta d\zzb\beta':\:\alpha,\alpha',
\beta,\beta'\in\Cinf\right\}
\]

\begin{prop}
$\left(\Omega^1,d\right)$ is a first order differential $*$-calculus over
$\Cinf$.
\end{prop}

The differential operators defined by \eqref{diffop} have some of the
properties of the classical $\frac{\partial}{\partial z}$ and
$\frac{\partial}{\partial\overline{z}}$ operators:

\begin{prop}
Let $f\in\Cinf$ be a restriction to $\Ybar$ of an entire function. Then
\[
\DRzb{f}=\DLzb{f}=0.
\]
Conversely, if $f\in\Cinf$ and either $\DRzb{f}=0$ or $\DLzb{f}=0$ then $f$
extends to an entire function.
\end{prop}

We also have an analogous result for antiholomorphic functions.

Let us recall the special function $\Fq$ introduced in \cite{e2}. It is a
continuous function on $\Ybar$ defined for
$z\not\in\left\{-q^{-2k}:\:k\in\ZZ_+\right\}$ by
\[
\Fq(z)=\prod_{k=0}^\infty\frac{1+q^{2k}\zbar}{1+q^{2k}z}.
\]
For reasons explained in \cite{opeq} (cf.~also \cite{nazb}) $\Fq$ is called the
{\em quantum exponential function.} This function is used in construction of
the multiplicative unitary of the quantum ``$az+b$'' group (\cite{azb})
as well as of the quantum $E(2)$ group (\cite{e2}). $\Fq$ does not belong to
$\Cinf$ because
the analytic extension of the map $\RR\ni t\mapsto\Fq(q^{it}z)$ has
singularities in the upper half plane (for generic $z$). Nevertheless we can
compute the action of the left differential operators on $\Fq$:
\[
\DLz{\Fq}=-\frac{1}{1-q^2}\Fq,\qquad
\DLzb{\Fq}=\frac{1}{1-q^2}\Fq.
\]
A slightly more sophisticated version of these formulas will be needed later:
\begin{equation}\label{diffFq}
\DLz{}\Fq(\zeta z)=-\frac{\zeta}{1-q^2}\Fq(\zeta z),\qquad
\DLzb{}\Fq(\zeta z)=\frac{\overline{\zeta}}{1-q^2}\Fq(\zeta z)
\end{equation}
for any parameter $\zeta\in\Ybar$.

It is also possible to construct higher order differential calculus on $\Ybar$.
The requirement that $d^2=0$ and that for any form $\rho$ the exterior
derivative $d\rho$ is given by the graded commutator with $\om-\omb$ leads to
the definition
\[
\Omega^2=\left.\left(\Omega^1\tens_{\Cinf}\Omega^1\right)\right/\mathscr{N},
\]
where $\mathscr{N}$ is the submodule generated by $d\zz\tens d\zz$, $d\zzb\tens
d\zzb$ and $d\zz\tens d\zzb+d\zzb\tens d\zz$. There are no higher dimensional
forms.

\subsection{Covariance}

Let us now relate the algebraic considerations of Subsection \ref{algdiff} and
the abstract definition of differential calculus of Subsection \ref{calcu}. as
before we shall identify the generator $\zz$ of the $\cst$-algebra $\CoYbar$
with the element $b$ affiliated with the $\cst$-algebra $A=\cpY$. Thus the
algebra generated by $b$ becomes a left coideal in the $\cst$-algebra with
comultiplication $\left(A,\DelG\right)$.

Now we need the functional $\ph$ defined on the algebra generated by $b$. We
cannot expect it to be continuous as it is designed to be the analogue of a
tangent vector. Let $\ph$ take the value $1$ on $b$ and $b^*$ and zero on all
higher powers of $b$ and $b^*$. Then $\ph$ satisfies \eqref{counit} in an
appropriate sense: the counit of the quantum group $G$ is a morphism from $A$ to
$\CC$ taking value $1$ on $a$ and $0$ on $b$. Since $\ph$ is equal to zero an
any product of more than one element we get an analogue of \eqref{counit}.

With this definition of $\ph$ we can formally apply $(\id\tens\ph)$ to
$\DelG(b)$ and call it $d\zz$:
\[
d\zz=\ph(b)a=a.
\]
Then the first formula of \eqref{dzom} gives
\[
\om=\frac{1}{q^2-1}b^{-1}a
\]
(notice that $b$ might not be invertible in some representations of $A$ and
consequently $b^{-1}a$ is not affiliated with $A$).

Now the rules of differential calculus can be given a precise analytical
meaning: for an element $f(b)$ of (the image in $A$ of) $\CoYbar$ the exterior
derivative $df$ is $\left[\om-\om^*,f\right]$. For functions $f$ belonging to
the algebra $\Cinf$ the rules of differentiating obtained in Subsection
\ref{calcu} are valid. Thus obtained differential calculus is covariant by
construction. However we are not claiming that a diagram analogous to
\eqref{HopfDiag} is commutative because the maps involved are not continuous.

The action of $G$ on differential forms on $\Hs=\Ybar$ is obtained by applying
$\DelG$ to the operator representing $d\zz$ affiliated with $A$. In other words
$\PhiG(d\zz)=a\tens d\zz$.

\section{Covariant measure}\label{covMeas}

Let us introduce a measure on $\Hs=\Ybar$. This measure will be covariant with
respect to the action of the quantum ``$az+b$'' group. For a positive
$f\in\CoYbar$ let
\[
\Int_\Ybar f\,
d\mu=\sum_{k\in\ZZ}\frac{q^{2k}}{2\pi}
\Int_0^{2\pi}f\left(q^ke^{i\ph}\right)\,d\ph.
\]

\begin{prop}\label{elem}
Let $f$ be a function in $\bB$ integrable with respect to $\mu$. Then we have
\[
\begin{array}{r@{\;=\;}l@{\medskip}}
\DS\Int_\Ybar f(qz)\,d\mu(z)&
\DS q^{-2}\Int_\Ybar f(z)\,d\mu(z),\\
\DS\Int_\Ybar f(q\cont z)\,d\mu(z)&
\DS\Int_\Ybar f(qz)\,d\mu(z).
\end{array}
\]
Inparticular the function $z\mapsto f(q\cont z)$ is integrable.
\end{prop}

\begin{cor}\label{norop}
Let $f$ and $g$ be functions in $\Cinf$. Then
\begin{enumerate}
\item\label{s1} if $f\DRz{g}$ and $\DLz{f}g$ are integrable then
\[
\Int_\Ybar f\DRz{g}\,d\mu=-q^2\Int_\Ybar \DLz{f}g\,d\mu;
\]
\item\label{s2} if $f\DRzb{g}$ and $\DLzb{f}g$ are integrable then
\[
\Int_\Ybar f\DRzb{g}\,d\mu=-q^{-2}\Int_\Ybar \DLzb{f}g\,d\mu.
\]

Let $\LYbar$ be the space of $\mu$-square integrable functions on $\Ybar$. The
scalar product in $\LYbar$ will be denoted by $\is{\cdot}{\cdot}$. With this
notation

\item\label{s3} if $f,g,\DRzb{f},\DRz{g}$ are in $\LYbar$ then
\[
\is{f}{\DRz{g}}=-q^2\is{\DRzb{f}}{g};
\]
\item\label{s4} if $f,g,\DLzb{f},\DLz{g}$ are in $\LYbar$ then
\[
\is{f}{\DLzb{g}}=-q^2\is{\DLz{f}}{g}.
\]
\end{enumerate}
\end{cor}

Let us define the following family of functions on $\Ybar$: for $k,l\in\ZZ$ a
nd $z\in\Ybar$ let
\begin{equation}\label{gkl}
\g{k,l}(z)=\left\{\begin{array}{c@{\quad\textrm{for}\quad}l@{\smallskip}}
\chi(z,q^l)&|z|=q^k,\\0&|z|\neq q^k.\end{array}\right.
\end{equation}
The family $\left(\g{k,l}\right)_{k,l\in\ZZ}$ is contained in $\Cinf$ and in
$\LYbar$. Moreover these functions form a maximal orthogonal system in $\LYbar$.
The action of the differential operators \eqref{diffop} on functions
$\g{k,l}$ can be computed:
\[
\begin{array}{r@{\;=\;}l@{\medskip}}
\DRz{\g{k,l}}&\DS\frac{1}{q^{-2}-1}\left(q^{-l-(k+1)}\;\g{k+1,l-1}
-q^{-k}\,\g{k,l-1}\right),\\
\DLz{\g{k,l}}&\DS\frac{1}{1-q^2}\left(q^{-k}\,\g{k,l-1}
-q^{l-(k-1)}\;\g{k-1,l-1}\right),\\
\DRzb{\g{k,l}}&\DS\frac{1}{1-q^2}\left(q^{-k}\,\g{k,l+1}
-q^{-l-(k-1)}\;\g{k-1,l+1}\right),\\
\DLzb{\g{k,l}}&\DS\frac{1}{q^{-2}-1}\left(q^{l-(k+1)}\;\g{k+1,l+1}
-q^{-k}\,\g{k,l+1}\right).
\end{array}
\]
We see that the linear span of functions $\left(\g{k,l}\right)_{k,l\in\ZZ}$ is
a dense subset of $\LYbar$ which is preserved by operators \eqref{diffop} and
left and right derivatives commute on this subset. Moreover, it follows from
\eqref{s3} and \eqref{s4} of Corollary \ref{norop} that operators \eqref{diffop}
are not only densely defined, but also closable. Without much trouble we obtain

\begin{prop}
Denoting by the same symbols the closures of the operators \eqref{diffop} in
the space $\LYbar$ we have
\[
\left(\DRz{}\right)^*=-q^2\DRzb{}\quad\text{and}\quad
\left(\DLz{}\right)^*=-q^{-2}\DLzb{}.
\]
Moreover $\DRz{},\DLz{},\DRzb{}$ and $\DLzb{}$ are normal operators on
$\LYbar$.
\end{prop}

There is another interesting result which is very much analogous to the
classical Stokes' theorem:

\begin{prop}
Let $f\in\Cinf$ have a compact support. Then
\[
\Int_{\Ybar}\DLz{f}\,d\mu=\Int_{\Ybar}\DRz{f}\,d\mu
=\Int_{\Ybar}\DLzb{f}\,d\mu=\Int_{\Ybar}\DRzb{f}\,d\mu=0.
\]
\end{prop}

Let us remark that $\mu$ is the only positive (non zero) measure on $\Ybar$ for
which the above theorem is true.

\subsection{Covariance}\label{covMu}

The covariance of $\mu$ can be expressed imprecisely as the property that
\begin{equation}\label{covar0}
(\id\tens\mu)\PhiG(f)=|a|^2\mu(f)
\end{equation}
for any $f\in\CoHs=\CoYbar$. For a more precise definition let us pass to the
von Neumann algebraic context. For this reason let us recall some basic facts
about the Haar measure on $G$.

The right Haar measure $h$ of $G$ was already introduced in \cite{vdH} (cf.~also
\cite{haar}). For an element $c=g(a)f(b)$ with $g\in\CoYhat$ and
$f\in\CoYbar$ we have
\[
h(c^*c)=\Int_{\Y}|g(\y)|^2\,d\y\Int_{\Ybar}|f(z)|^2\,d\mu(z),
\]
where $d\y$ denotes the Haar measure on $\Y$:
\begin{equation}\label{HaarY}
\Int_{\Y}g(\y)\,d\y=\sum_{k\in\ZZ}
\frac{1}{2\pi}\Int_0^{2\pi}g\left(q^ke^{it}\right)\,dt
\end{equation}
(notice that $\mu$ is the push forward of the measure $|\y|^2d\y$ via the
inclusion of $\Y$ into $\Ybar$).

Let us represent $A$ in the GNS Hilbert $H$ space for the weight $h$
(one can show that $H=\LY\tens\LYbar$ where $\LY$ is the space of functions on
$\Y$ square integrable with respect to the Haar measure \eqref{HaarY}). We can
form the von Neumann algebra $M=A''$. Then $h$ extends to a normal faithful
semifinite weight on $M$ (cf.~\cite{lcqg} or \cite{mnw}). One can show that
$M$ is the von Neumann algebra crossed product of $N=\LinfYbar$ (where on
$\Ybar$ we take the measure class of $\mu$) by the natural action of $\Y$,
which clearly extends the action $\beta$ on $\CoYbar$.
The weight $h$ is invariant under the dual action of $\Yhat$ on $M$. Indeed, the
Haar measure $h$ is right invariant, i.e.~for any $\ph\in A^*_+$ and any
$x\in A_+$ we have $h(\ph*x)=\ph(I)h(x)$. In the proof of Theorem \ref{HomSpace}
we saw that the dual action of $\Y$ on $A$ is given by convolution with
functionals $\pi_{\yhat}$. Therefore denoting the dual action by $\behat$ we
have $h=h\comp\behat_{\yhat}$ for any $\yhat\in\Yhat$ (in the proof of Theorem
\ref{HomSpace} we denoted $\behat$ by $\alhat$). It follows that the extensions
$h$ and $h\comp\behat_{\yhat}$ to $M$ are equal for all $\yhat\in\Yhat$. As the
dual action on $M$ is the natural extension of that on $A$ we have that the
normal extension of $h$ is invariant under the (von Neumann algebraic) dual
action.

Let $T$ be th operator valued weight $M\to N$ coming from the integration of
the dual action (\cite{HaaOp}). Then by \cite[Theorem 3.7]{HaaDuI} and
\cite[Theorem 1.1]{HaaDuII} There exists a unique measure $\nu$ on $\Ybar$ such
that $h=\nu\comp T$. One easily checks that for $c=g(a)f(b)\in A\subset M$
we have $\bigl(\nu\comp T\bigr)(c^*c)=\bigl(\mu\comp T\bigr)(c^*c)$. By
uniqueness of $\nu$ we have $\nu=\mu$ and $h=\mu\comp T$.

Let $\rho$ be the modualr element of $G$ (the Radon-Nikodym derivative of the
left Haar measure $h^{\text{\rm\tiny L}}=h\comp R$ with respect to $h$). Then
$\rho$ is a positive selfadjoint element affiliated with $A$, but we can also
treat it as an unbounded operator on $H$
By \cite[Proposition 2.5]{unit} we have
\begin{equation}\label{covar}
\left((\omega_{\xi,\xi}\tens\id)\PhiG(f)\right)=\|\rho\xi\|^2\mu(f)
\end{equation}
for all $f\in\CoHs$ and all $\xi\in D(\rho)$.

Since $\rho=|a|^2$ we shall take \eqref{covar} as the precise formulation of
\eqref{covar0}.

\section{Fourier transform}\label{FTs}

Our definition of the Fourier transformation on $\Ybar$ will be motivated by
the properties of the quantum exponential function $\Fq$. Let $f$ be an element
of $\Cinf$ which is integrable for the measure $\mu$. Then we define
\[
\left(\F f\right)(\zeta)=\Int_\Ybar\Fq(\zeta z)f(z)\,d\mu(z).
\]
In order to describe the properties of the operation $\F$ let us introduce the
following notation: for $t\in\RR$ and any function $g$ on $\Ybar$ let
\[
\begin{array}{r@{\;=\;}l@{\smallskip}}
\left(\sigma_t(g)\right)(z)&g\left(q^{-it}z\right),\\
\left(\theta(g)\right)(z)&g\left(qz\right).
\end{array}
\]
Then $\left(\sigma_t\right)_{t\in\RR}$ becomes a one parameter group of
automorphisms of the algebra of functions on $\Ybar$ and similarly $\theta$ is
an automorphism of this algebra commuting with the automorphisms
$\left(\sigma_t\right)_{t\in\RR}$. Restricting attention to the
continuous functions the group $\left(\sigma_t\right)_{t\in\RR}$ becomes
continuous in the natural topology and we can form the analytic generator
$\sigma_i$. The domain of $\sigma_i$ is precisely the algebra $\bB$ described in
Subsection \ref{smooth}. From Proposition \ref{elem} it easily follows that
\begin{equation}\label{male}
\begin{array}{r@{\;=\;}l@{\smallskip}}
\left(\sigma_{-i}\comp\theta^{-1}\right)\comp\F&
q^2\F\comp\left(\sigma_i\comp\theta\right),\\
\left(\sigma_{-i}\comp\theta\right)\comp\F&
q^{-2}\F\comp\left(\sigma_i\comp\theta^{-1}\right).
\end{array}
\end{equation}
Indeed, the fact that $\theta\comp\F=q^{-2}\F\comp\theta$ is a direct
consequence of the first formula of Proposition \ref{elem}. In order to see
that $\sigma_{-i}\comp\F=\F\comp\sigma_i$ let us compute
\[
\begin{array}{r@{\;=\;}l@{\medskip}}
\left(\F f\right)(q^{-it}\zeta)&
\DS\Int_\Ybar\Fq\left(q^{-it}\zeta z\right)f(z)\,d\mu(z)\\
&\DS\Int_\Ybar\Fq\left(\zeta z\right)f(q^{it}z)\,d\mu(z)
\end{array}
\]
Now if $f$ and $\sigma_i(f)$ are integrable then the last expression above has
a holomorphic continuation to $t=-i$. Therfore
$\left(\F f\right)(q^{-1}\cont\zeta)$ exists and is equal to
$\left(\F\comp\sigma_i(f)\right)(\zeta)$.

Let us denote the operator of multiplication by the identity function $\zz$ on
$\Ybar$ by the symbol $\zZ$ and let $\zZ^*$ be the multiplication by $\zzb$.

\begin{thm}
The commutation relations between $\F$, $\zZ$, $\zZ^*$ and the operators
\eqref{diffop} are the following:
\[
\begin{array}{r@{\;=\;}l@{\qquad}r@{\;=\;}l@{\medskip}}
\zZ\comp\F&(q^{-2}-1)\,\F\comp\DRz{},&\zZ^*\comp\F&-(q^2-q^4)\,\F\comp\DRzb{},
\\
\DLz{}\comp\F&-\DS\frac{1}{1-q^2}\,\F\comp\zZ,&
\DLzb{}\comp\F&\DS\frac{1}{1-q^2}\,\F\comp\zZ^*,\\
\DRz{}\comp\F&
-\DS\frac{q^4}{1-q^2}\,
\F\comp\zZ\comp\left(\sigma_i\comp\theta\right),
&\DRzb{}\comp\F&\DS\frac{q^4}{1-q^2}\,
\F\comp\zZ^*\comp\left(\sigma_i\comp\theta^{-1}\right).
\end{array}
\]
\end{thm}

\begin{proof}
The formulas in the first line are easy consequences of \eqref{diffFq} and
statements \eqref{s1} and \eqref{s2} of Proposition \ref{norop}. We shall only
give the proof of the first one:
\[
\begin{array}{r@{\;=\;}l@{\medskip}}
\left(\zZ\comp\F(f)\right)(\zeta)
&\zeta\DS\Int_\Ybar\Fq(\zeta z)f(z)\,d\mu(z)
=\DS\Int_\Ybar\zeta\Fq(\zeta z)f(z)\,d\mu(z)\\
&-(1-q^2)\DS\Int_\Ybar\frac{-\zeta}{1-q^2}
\Fq(\zeta z)f(z)\,d\mu(z)\\
&-(1-q^2)\DS\Int_\Ybar\DLz{}\Fq(\zeta z)f(z)\,d\mu(z)\\
&(1-q^2)q^{-2}\DS\Int_\Ybar\Fq(\zeta z)\DRz{f}(z)\,d\mu(z)
=(q^{-2}-1)\left(\F\DRz{f}\right)(\zeta).
\end{array}
\]
The two formulas in the last line are obtained from the two in the second line
by using relations \eqref{LR} and \eqref{male}. Finally the formulas in the
second line are obtained from \eqref{diffFq}. Again we only give proof of
the first one:
\[
\begin{array}{r@{\;=\;}l@{\medskip}}
\left(\DS
\frac{\partial_{\scriptscriptstyle\mathrm{L}}}{\partial\zeta}\F f\right)(\zeta)
&\DS\frac{\partial_{\scriptscriptstyle\mathrm{L}}}{\partial\zeta}
\Int_\Ybar\Fq(\zeta z)f(z)\,d\mu(z)\\
&\DS\Int_\Ybar\frac{\partial_{\scriptscriptstyle\mathrm{L}}}{\partial\zeta}
\Fq(\zeta z)f(z)\,d\mu(z)\\
&\DS\Int_\Ybar\frac{-\zeta}{1-q^2}\Fq(\zeta z)f(z)\,d\mu(z)\\
&\DS\frac{-\zeta}{1-q^2}\Int_\Ybar\Fq(\zeta z)f(z)\,d\mu(z)
=-\frac{1}{1-q^2}\left(\F\zZ f\right)(\zeta).
\end{array}
\]
\end{proof}

It is not difficult to show that the operator $\F$ can be extended to a
closed operator on $\LYbar$. The adjoint $\F^*$ is also an integral operator
with kernel $\overline{\Fq}$. With a detailed analysis of the special function
$\Fq$ one can prove the following theorem:

\begin{thm}
\noindent
\begin{enumerate}
\item The functions $\left(\F\g{k,l}\right)_{k,l\in\ZZ}$ introduced by
\eqref{gkl} are in $\Cinf$ and are integrable with respect to $\mu$.
\item We have $\F^*\F\g{k,l}=(q\sigma_{-i})^2\g{k,l}$.
\item The operator $q^{-1}\F\comp\sigma_{-i}$ extends to a unitary operator on
$\LYbar$.
\end{enumerate}
\end{thm}

Notice that if we formally put $q=1$ then the group
$\left(\sigma_t\right)_{t\in\RR}$ is trivial and consequently $\F$ becomes
a unitary operator.

\end{document}